\documentclass[10pt]{article}
\font\smallit=cmti10
\font\smalltt=cmtt10

\usepackage{comment}
\usepackage{amssymb,latexsym,amsmath,epsfig,amsthm} 
\usepackage{hyperref}
\makeatletter

\renewcommand\section{\@startsection {section}{1}{\z@}
{-30pt \@plus -1ex \@minus -.2ex}
{2.3ex \@plus.2ex}
{\normalfont\normalsize\bfseries\boldmath}}

\renewcommand\subsection{\@startsection{subsection}{2}{\z@}
{-3.25ex\@plus -1ex \@minus -.2ex}
{1.5ex \@plus .2ex}
{\normalfont\normalsize\bfseries\boldmath}}

\renewcommand{\@seccntformat}[1]{\csname the#1\endcsname. }

\newcommand\blfootnote[1]{
\begingroup
\renewcommand\thefootnote{}\footnote{#1}
\addtocounter{footnote}{-1}
\endgroup
}

\makeatother

\newtheorem{theorem}{Theorem}

\newtheorem{proposition}{Proposition}
\newtheorem{corollary}{Corollary}

\theoremstyle{definition}

\newtheorem{conjecture}{Conjecture}

\begin{document}														\begin{figure}
\vspace*{-62pt}
\vspace{.75in}
\end{figure}
\vspace*{-45pt}
\leftline{\smalltt\#A1} \vskip -12.5pt
\centerline{\smalltt  INTEGERS 26 (2026)}
\vskip 3pt \hrule												

\begin{center}
\uppercase{\bf \boldmath On the Distribution of Points of Valuation $1$ for a Polynomial in Two Variables}\blfootnote{DOI: }
\vskip 20pt
{\bf Krishnan Rajkumar}\\
{\smallit School of Engineering, Jawaharlal Nehru University, New Delhi, India }\\
{\tt krishnan@jnu.ac.in}\\ 
\vskip 10pt
{\bf Shubham}\\
{\smallit School of Physical Sciences, Jawaharlal Nehru University, New Delhi, India}\\
{\tt shubha76\_sps@jnu.ac.in}\\ {\tt shubham01nitw@gmail.com}
\end{center}
\vskip 20pt
\centerline{\smallit Received: , Revised: , Accepted: , Published: } 
\vskip 30pt

\centerline{\bf Abstract}
\noindent
We investigate the variation in the total number of points in a random $p\times p$ square in $\mathbb{Z}^2$ at which the $p$-adic valuation of a given polynomial in two variables is greater than $1$. We establish that  this quantity follows a Poisson distribution as $p\rightarrow\infty$ under a certain conjecture. We also relate this conjecture to certain uniform distribution properties of a vector-valued sequence.

\pagestyle{myheadings}
\markright{\smalltt INTEGERS: 26 (2026)\hfill}
\thispagestyle{empty}
\baselineskip=12.875pt
\vskip 30pt

\section{Introduction}\label{sec1}
The $p$-adic valuation of integer sequences is a fundamental concept in number theory, providing insights into the divisibility properties of integers with respect to a prime number \( p \). The {\it $p$-adic valuation} $\nu_p(n)$ of an integer $n$ is defined as the largest exponent $e$ such that $p^e$ divides $n$ for a prime number $p$. 
Historically, the concept of $p$-adic valuation stems from the work of Hensel \cite{four}, who introduced $p$-adic numbers in 1897. This development allowed mathematicians to study numbers and their divisibility in a local sense, focusing on individual primes. Over time, the theory of $p$-adic numbers and valuations was expanded by E. Artin, Grothendieck, and others, making it a crucial tool in areas such as Diophantine equations and algebraic geometry.

One of the first and important results in the study of $p$-adic valuations is Legendre’s formula \cite{five} for the $p$-adic valuation of \( n! \). For a prime number \( p \), the formula is given by:

\[
\nu_p(n!) = \sum_{k=1}^{\infty} \left\lfloor \frac{n}{p^k} \right\rfloor
,\]where $\lfloor x\rfloor$ is the greatest integer less than or equal to $x$. This sums the number of multiples of \( p, p^2, p^3, \dots \) less than or equal to \( n \), counting the number of times \( p \) divides \( n! \). This formula has important applications in number theory and combinatorics, particularly in understanding binomial coefficients and their divisibility by primes.

The work presented by Byrnes et al. \cite{two} studies the $2$-adic valuation of $n^2-a$.  It is shown that $\nu_2(n^2-a)$ has a simple closed form (see \cite{one,three} for a detailed discussion of the definition of a closed form) when $a \not \equiv 4,7 \mod 8$.  For these two remaining cases, the valuation cannot be described in a simple form. However, it can be fully described in a rule-based framework called the \textit{valuation tree}. 

These examples form part of a general project, with several notable contributions made by Victor Moll~\cite{two, moll} and his coauthors to describe the $p$-adic valuations $\nu_p(x_n)$, 
for any given sequence $\{x_n\}_{n \in \mathbb{N}} $ of integers. In the work \cite{seven} by the second author, the valuations $\nu_p(f(m,n))$ of a two-dimensional sequence $\{f(m,n)\}_{m,n \in \mathbb{N}}$ of integers are taken up, 
where $f(x,y) \in \mathbb{Z}[x,y]$ is any polynomial in two variables with integer coefficients. It is fully described by generalizing the notion of the valuation tree. 

In this work, we wish to understand the distribution of the $p$-adic valuations $\nu_p(f(m_1,m_2,\ldots,m_k))$ as $m_1,m_2,\ldots,m_k$ vary over $\mathbb{N}$ for a polynomial $f \in \mathbb{Z}[x_1,x_2,$ $\ldots,$ $ x_k]$. Let us first discuss the case $k=1$ of a polynomial in one variable.

Let $f(x)\in \mathbb{Z}[x]$ have degree $d$ and non-zero discriminant. First, we note that points $n$ where $f$ has non-zero valuation $\nu_p(f(n))>0$ follow a periodic pattern with period $p$, as $\nu_p(f(x))>0 \Leftrightarrow \nu_p(f(x+p))>0$. Hence, every interval of length $p$, in particular the intervals $\{kp,kp+1,\ldots,kp+p-1\}$, will contain the same number of points of non-zero valuation and this number is the cardinality of the solution set of $f(x)=0$ over $\mathbb{F}_p$. 

The situation changes when we consider the points of valuation $\nu_p(f(n))>1$. In the first place, the pattern is still periodic with period $p^2$. However, not all intervals of length $p$ contain the same number of points where $\nu_p(f(n))>1$. In particular, at least $p-d$ intervals $\{kp,kp+1,\ldots,kp+p-1\}$ with $k=0,1,\ldots,p-1$ do not contain any points where $\nu_p(f(n)) > 1$. This is because, by Corollary \ref{cor11}, the total number of points where $\nu_p(f(n))>1$ in the set $\{0,1,\ldots,p^2-1\}$ is at most $d$. 
Hence, as $p\rightarrow\infty$, the total number of points $n$ of valuation $\nu_p(f(n))>1$ in an interval $\{kp,kp+1,\ldots,kp+p-1\}$ chosen uniformly at random from $k=0,1,\ldots,p-1$ is zero with probability tending to $1$.

We now turn to the case $k=2$ of a polynomial in two variables and answer the following question: what is the distribution of points in $\mathbb{N}^2$ at which a polynomial $f(x,y) \in \mathbb{Z}[x,y]$ has valuation greater than $1$?

For a fixed polynomial $f(x,y) \in \mathbb{Z}[x,y]$, we note that the points $(r,s)$ where $f$ has non-zero valuation $\nu_p(f(r,s))>0$ display a $p$-periodic pattern in both $r$ and $s$ coordinates. The number of points in any $p \times p$ square is the number of points on the curve $C_f(\mathbb{F}_p) = \{(x,y) \in \mathbb{Z}^2: \nu_p(f(x,y)) > 0 \text{ and } 0\leq x,y <p\}$ given by the equation $f(x,y)=0$ over $\mathbb{F}_p$. Let $m=|C_f(\mathbb{F}_p)|$ be the number of points on this curve.

\begin{figure}[h]  
    \centering  
   \includegraphics[width=0.8\textwidth]{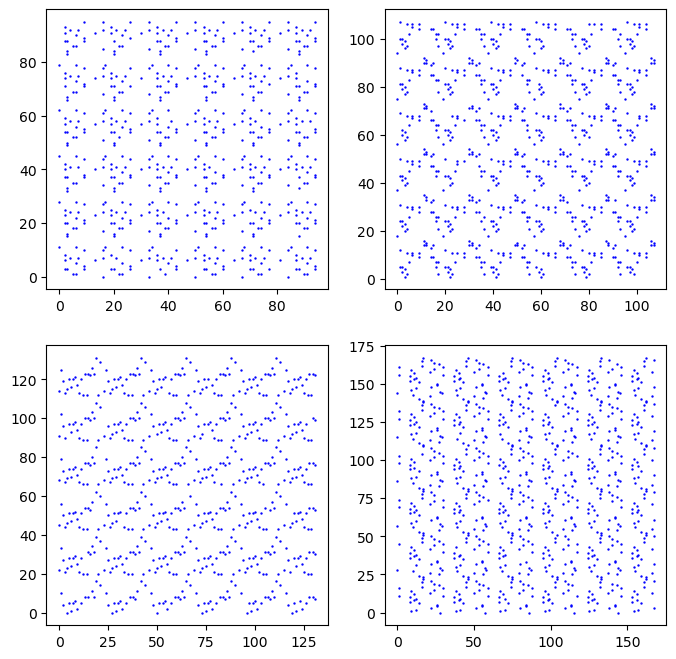}  
   \caption{Scatter Plot of $(x, y)$ satisfying: \textbf{Upper Left}: $\nu_p(x^3 + y^3 + x^2 y + y + 1) = 1$, for $p = 17$; \textbf{Upper Right}: $\nu_p(x^3 + y^2 x + xy + x + y + 1) = 1$, for $p = 19$; \textbf{Lower Left}: $\nu_p(x^3 + xy + x + y + 1) = 1$, for $p = 23$; \textbf{Lower Right}: $\nu_p(y^2 x + xy + x + y + 1) = 1$, for $p = 29$.}
    \label{fig:1}  
\end{figure}

Let us observe the patterns that occur for the points in $\mathbb{N}^2$ at which the valuation of $f$ is exactly $1$, for various polynomials $f$. This is shown in Figure \ref{fig:1}.
The points missing from $C_f(\mathbb{F}_p)$ in the various $p \times p$ squares seem to exhibit some element of randomness. In order to quantify this, we first study the properties of indicator random variables $Y_{x,y}$ in Section \ref{sec3}. In summary, for each $(x,y) \in C_f(\mathbb{F}_p)$, we define a random variable $Y_{x,y}$ that indicates which translates of $(x,y)$ into a random $p \times p$ square have valuation greater than $1$. These random variables are identically distributed Bernoulli random variables with probability $1/p$, though they are not always pairwise independent.

In order to further quantify this variation, we define the following random variable $X : \mathbb{F}_p \times \mathbb{F}_p \rightarrow \mathbb{Z}$:
\begin{align} \label{Xdefn}
X(k,l) 
&= \bigl|\bigl\{(x,y) \in 
   \{kp, kp+1, \ldots, kp+p-1\} \times \{lp, lp+1, \ldots, lp+p-1\} \nonumber\\
&\qquad\qquad : \ \nu_p(f(x,y)) > 1 \bigr\}\bigr|.
\end{align}

\begin{figure}[h]  
    \centering  
    \includegraphics[width=0.8\textwidth]{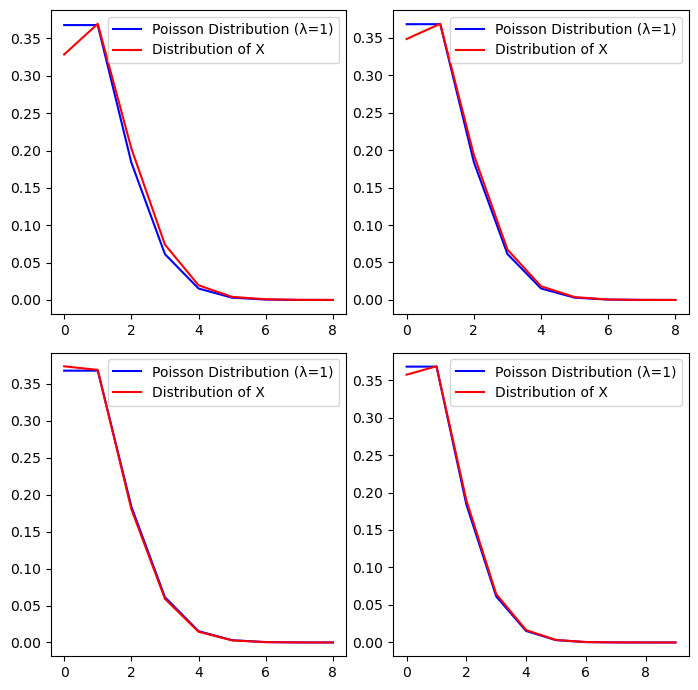}  
   \caption{Comparison of the Poisson distribution ($\lambda=1$) with the distribution of $X$ for $f(x,y)=x^3+y^2+x y+1$ and the prime numbers, \textbf{Upper Left}: $p = 211$; \textbf{Upper Right}: $p = 503$; \textbf{Lower Left}: $p = 1013$; \textbf{Lower Right}: $p = 1033$.}
    \label{fig:2}  
\end{figure}

\noindent In other words, $X$ represents the total number of points in a random $p \times p$ square at which the valuation of $f$ is greater than $1$. The distribution of $X$ for various primes $p$ is shown in Figure \ref{fig:2}.

The distribution of $X$ appears to approach the Poisson distribution for large $p$. This limiting behaviour of $X$ can be explained by the following reasoning under the assumption that $Y_{x,y}$ are mutually independent: by \eqref{X_def} and this assumption, we get that $X$ follows the binomial distribution with parameters $m$ and $1/p$, with their product $m/p$ tending to $1$ as $p\rightarrow \infty$ (by the Hasse-Weil bounds \cite{eight}). Hence, under this assumption, the distribution of $X$ tends to a Poisson distribution with the parameter $\lambda=1$. The central result of this work is that this statement holds without the assumption on mutual independence of $Y_{x,y}$, under Conjecture \ref{conj1}. To state this conjecture, we need some notation.

For any $k$-tuple $\mathbf{P}=((r_1,s_1),(r_2,s_2),\dots, (r_k,s_k))$ of points on the curve $C_f(\mathbb{F}_p)$, denote the associated matrix $B_k(\mathbf{P})$ with entries in $\mathbb{F}_p$ as
\begin{equation}\label{B_k}
	B_k(\mathbf{P})= \begin{bmatrix}
		f_x(r_1,s_1) &  f_y(r_1,s_1)& f(r_1, s_1)/p \\
		f_x(r_2,s_2) & f_y(r_2,s_2) & f(r_2, s_2)/p\\
		\vdots & \vdots & \vdots \\
		f_x(r_k,s_k) & f_y(r_k,s_k) & f(r_k, s_k)/p 
	\end{bmatrix}.
\end{equation}

For example, consider $f(x,y)=x^3+y^2+xy+1$ and $p=7$. Then the curve $C_f(\mathbb{F}_p)=\{(1, 3), (3, 0), (3, 4), (4, 1), (4, 2), (5, 0), (5, 2), (6, 0), (6, 1)\}$.
\noindent
For $k=3$, and $\mathbf{P_1}= ((1,3), (3,0), (3,4))$, we have
\begin{equation*}
	B_3(\mathbf{P_1})= \begin{bmatrix}
		6 &  7 & 2 \\
		27 & 3 & 4\\
		31 & 11 & 8 
	\end{bmatrix}=\begin{bmatrix}
	  6 & 0 & 2 \\
    6 & 3 & 4 \\
    3 & 4 & 1
	\end{bmatrix},
\end{equation*}
in $\mathbb{F}_7$, which has rank $3$ (its determinant is non-zero in $\mathbb{F}_7$).

For $\mathbf{P_2}= ((1,3), (3,0), (5,2))$, we have
\begin{equation*}
	B_3(\mathbf{P_2})= \begin{bmatrix}
		6 &  7 & 2 \\
		27 & 3 & 4\\
		77 & 9 & 20 
	\end{bmatrix} = \begin{bmatrix}
	  6 & 0 & 2 \\
    6 & 3 & 4 \\
    0 & 2 & 6
	\end{bmatrix},
\end{equation*}
ov $\mathbb{F}_7$, which has rank $2$, as the first two columns are linearly independent (any of the $2 \times 2$ minors involving the first two columns is non-zero in $\mathbb{F}_7$) and the third column is a linear combination of the first two columns over $\mathbb{F}_7$.  

Let $m_{k}(f,p)$ denote the number of $k$-tuples $\mathbf{P} \in C_f(\mathbb{F}_p)^k$ that are element-wise distinct, for which $B_k(\mathbf{P})$ has rank $2$ and the first two columns of $B_k(\mathbf{P})$ are linearly independent. In the preceding example, $m_3(f,p)$ includes $\mathbf{P_2}$ and its permutations, but not $\mathbf{P_1}$ or its permutations.

For $k=1$, we observe that $m_{k}(f,p)=0$ trivially. 
For $k\geq 2$, 
 we have the following conjecture. 
\begin{conjecture}\label{conj1} For all $k\geq 2$, we have
	\begin{align}
		&\lim_{p\to\infty}  \frac{m_{k}(f,p)}{p^2} = 1. \label{lim1}
	\end{align}
\end{conjecture}

\begin{figure}[h]  
    \centering  
   \includegraphics[width=0.8\textwidth]{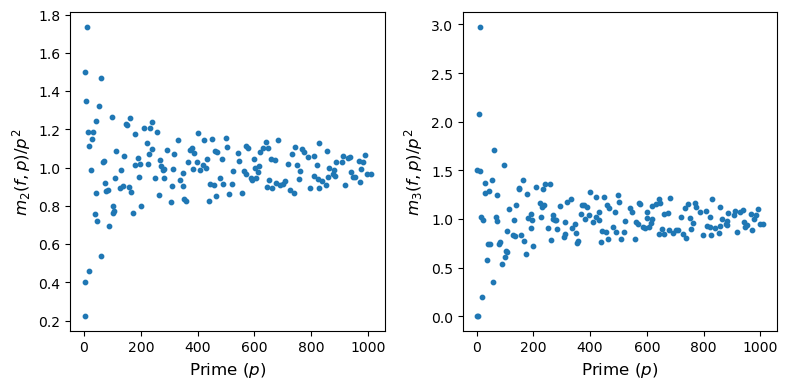}  
   \caption{\textbf{Left}: Plot of $m_{2}(f,p)/p^2$ versus $p$; \textbf{Right}: Plot of $m_{3}(f,p)/p^2$ versus $p$ for the polynomial $f(x,y)=x^3+y^2+x y+1$. }

    \label{fig:3}  
\end{figure}

Qualitatively, this conjecture implies that, for $k>2$, the $k$-tuples counted by $m_k(f,p)$ are very sparse among the element-wise distinct $k$-tuples in $C_f(\mathbb{F}_p)^k$. In continuation of the preceding example with $f(x,y)=x^3+y^2+xy+1$ and $p=7$, we calculate $m_3(f,7)=90$, which is comparable to $p^2=49$, although their ratio is much higher than $1$. However, for larger values of $p$, numerical observations appear to corroborate this conjecture, as seen in Figure \ref{fig:3}.
 We also show in Section \ref{sec4} that this conjecture is a natural consequence of certain versions of the uniform distribution of the sequence $\mathcal{F}_p$ in $\mathbb{F}_p^3$, where
\begin{equation}\label{Fcal_defn}
	\mathcal{F}_p = \{ (f_x(r,s),f_y(r,s),f(r,s)/p) \in \mathbb{F}_p^3 \ | \ (r,s) \in C_f(\mathbb{F}_p)\}.
\end{equation}

Now we are in a position to state the main result of this work.

\begin{theorem}\label{thm1}
	Let $f(x,y)\in\mathbb{Z}[x,y]$ be a polynomial with no solutions to the set of equations $f(x,y)\equiv 0 \pmod p, f_x(x,y)\equiv 0 \pmod p, \text{ and } f_y(x,y)\equiv 0 \pmod p$ for all sufficiently large prime numbers $p$. Let $X$ be the random variable defined in \eqref{Xdefn}, and assume that Conjecture \ref{conj1} is satisfied. Then X converges in distribution to a Poisson random variable with parameter $\lambda=1$ as $p\rightarrow\infty$.
\end{theorem}

We note that the condition on the polynomial $f(x,y)$ mentioned in Theorem \ref{thm1} is satisfied by the polynomials mentioned in Figures \ref{fig:1}, \ref{fig:2}, and \ref{fig:3}. For instance, the polynomial $f(x,y)=x^3+y^2+x y+1$ does not have any solutions to the set of equations $f(x,y)\equiv 0 \pmod p,\ f_x(x,y)\equiv 0 \pmod p, \text{ and } f_y(x,y)\equiv 0 \pmod p$ for $p>431$.

 To obtain a sufficient condition on $f(x,y)$ that is independent of $p$, one can eliminate the variable $y$ from $f(x,y)= 0, \ f_x(x,y)= 0$, as well as from $f(x,y)= 0, \ f_y(x,y)= 0$. This would result in two polynomials in $x$, both of which have roots modulo $p$ at the solutions of the original set of equations. However, there can be common roots for these two polynomials modulo $p$ if and only if $p \mid R$, where $R$ is the resultant of these two polynomials. Thus, a sufficient condition on $f(x,y)$ for Theorem \ref{thm1} to hold is that this resultant $R$ is a non-zero integer.
 
 We also remark that the case of polynomials in more than two variables is considerably involved and will be considered in future research.

We now describe the organization of this paper. In Section \ref{sec2}, we prove several elementary results on the number of solutions to certain valuation conditions. In Section \ref{sec3}, we interpret these results as describing the single and joint distributions of the indicator random variables indexed by points in $C_f(\mathbb{F}_p)$. We also derive the implication of Conjecture \ref{conj1} on the sum of expectations of $k$-fold products of these random variables. Section \ref{sec4} contains results on the equidistribution of $\mathcal{F}_p$ and their  relation to Conjecture \ref{conj1}. Section \ref{sec5} contains the proof of Theorem \ref{thm1}. 
 
\section{Preliminaries}\label{sec2}
We start by studying the cardinality of the set $\{x : \nu_p(g(x))>1, 0\leq x < p^2\} $ for a given polynomial $g(x)\in \mathbb{Z}[x]$.
\begin{proposition}\label{Pre 0}
	Let $p$ be a prime number, let $g(x) \in \mathbb{Z}[x]$, and let $r \in \{0,1,\ldots,p-1\}$ not be a common solution to equations $g(r)\equiv 0 \pmod p \text{ and } g'(r) \equiv 0 \pmod p$. Then the cardinality of the set $ \{ k \in \mathbb{F}_p : \nu_p(g(r+kp)) >1\}$ is at most $1$.
\end{proposition}
\begin{proof}
	We note that \begin{equation}\label{1var}
		g(r+kp) \equiv g(r)+kpg'(r) \pmod {p^2}. 
	\end{equation}
    
	If $g(r) \not\equiv 0 \pmod p$, then by \eqref{1var}, there are no solutions to $g(r+kp) \equiv 0 \pmod {p^2}$.   
	On the other hand, if $g(r) \equiv 0 \pmod p$ and $g'(r) \not\equiv 0 \pmod p$, then \eqref{1var} provides a unique solution to $g(r+kp)\equiv 0 \pmod {p^2}$,  where $k\in \mathbb{F}_p$.	
\end{proof}

\begin{corollary}\label{cor11}
Let $g(x)\in \mathbb{Z}[x]$ have degree $d$ and non-zero discriminant $\Delta$. Then the cardinality of the set $S = \{ x : \nu_p(g(x)) > 1,\ 0 \le x < p^2 \}
$ is at most $d$ for sufficiently large prime numbers $p$.
\end{corollary}
\begin{proof}
	First, we note that for $p$ sufficiently large, we have $\Delta \not\equiv 0 \pmod p$, as $\Delta$ is a non-zero integer (since it is a polynomial in the coefficients of $g$). Hence, there are no repeated roots of $g$ in any algebraically closed extension of $\mathbb{F}_p$. As a result, there are no solutions to $g(r)\equiv 0 \pmod p$ and $g'(r) \equiv 0 \pmod p$ in $\mathbb{F}_p$. 
	
	Thus, the conditions of Proposition \ref{Pre 0} hold for all $r\in \mathbb{F}_p$. Now, from the proof of Proposition \ref{Pre 0}, it is clear that there are solutions to $g(r+kp) \equiv 0 \pmod {p^2}$ only if $r$ is a root of $g(r)\equiv 0 \pmod p$. Hence, there are at most $d$ such roots $r \in \{0,1,\ldots,p-1\}$ to consider, and at most one translate $r+kp$ belongs to the set $S$. Thus, the cardinality of $S$ is at most $d$.
\end{proof}

Next, we consider a polynomial $f(x,y)\in\mathbb{Z}[x,y]$ and study the cardinality of the set $\{(x,y) : \nu_p(f(x,y)) = 1,\ 0 \le x,y < p^2\}$. We get the following results.
\begin{proposition}\label{Pre 1}
Let $p$ be a prime number, let $f(x,y) \in \mathbb{Z}[x,y]$ be a polynomial and let $(r,s)\in C_f(\mathbb{F}_p)$ such that at least one of $f_x(r,s)$ or $f_y(r,s)\not\equiv 0 \pmod p$. Then the cardinality of the set $\{(k,l) \in \mathbb{F}_p^2 : \nu_p(f(r+kp, s+lp)) = 1\}$ is $p(p-1)$.
\end{proposition}
\begin{proof}
Consider the solutions of \begin{equation*} 
  f(r+kp, s+lp)   \equiv 0 \pmod{p^2} ,
\end{equation*}for $k,l \in \{0,1,2,\ldots,p-1\}$. Now by Taylor series expansion, we have \begin{equation}
f(r,s)+kpf_x(r,s)+lpf_y(r,s) \equiv 0 \pmod {p^2}.  \label{1}
\end{equation} We know that $f(r,s) \equiv 0 \pmod p$ which implies there exists an integer $\alpha$, such that $f(r,s) = \alpha p$. Now, by using \eqref{1}, we have \begin{equation}\label{2}
\alpha+kf_x(r,s)+lf_y(r,s) \equiv 0 \pmod p.
\end{equation}Assuming, without loss of generality, that $f_x(r,s)\not\equiv 0 \pmod p$, we have \begin{equation*}
k \equiv \frac{-\alpha-lf_y(r,s)}{f_x(r,s)} \pmod p.
\end{equation*} So, there are $p$ choices for $(k,l)$ such that $f(r+kp, s+lp)   \equiv 0 \pmod{p^2} $. Hence, there are $p^2-p$ choices for $(k,l)$ such that $f(r+kp, s+lp)  \not \equiv 0 \pmod{p^2} $.
\end{proof}

\begin{corollary}
Let $m$ and $s$ be the cardinalities of the sets $C_f(\mathbb{F}_p)$ and $S$ respectively for a prime number $p$, where $S = \{(x,y): \nu_p(f(x,y))=1, 0\leq x,y < p^2\} $. Then $s =m p(p-1)$.
\end{corollary}
We also derive the following result directly from \eqref{2}.
\begin{corollary}\label{cor3}
Let $p$ be a prime number and let $(r,s)\in C_f(\mathbb{F}_p)$ be such that $f_x(r,s) \equiv f_y(r,s) \equiv 0 \pmod p$. Then the cardinality of the set
$\{ (k,l) \in \mathbb{F}_p^2 : \nu_p(f(r+kp, s+lp)) = 1 \}$
is $p^2$ if $f(r,s) \equiv 0 \pmod {p^2}$, and $0$ if $f(r,s) \not\equiv 0 \pmod {p^2}$.
\end{corollary}

Next, we examine the mutual influence of points $(r_1,s_1)$ and $(r_2,s_2)$ on the curve $C_f(\mathbb{F}_p)$.  
 Recall the matrix $B_2((r_1,s_1),(r_2,s_2))$ from \eqref{B_k}, denoted as $A$ for convenience, 
\begin{equation*}
A = \begin{pmatrix}
	f_x(r_1,s_1) &  f_y(r_1,s_1)& f(r_1,s_1)/p \\
	f_x(r_2,s_2) & f_y(r_2,s_2) & f(r_2,s_2)/p
\end{pmatrix}\cdot
\end{equation*} Let $r$ be the rank of $A$. In the case where $r=2$ and the first two columns of $A$ are linearly independent, we establish the following.
 \begin{proposition}\label{prop3}
 Let $p$ be a prime number, let $f(x,y) \in \mathbb{Z}[x,y]$ be an integer polynomial and let $(r_1,s_1)$ and $(r_2,s_2)$ be points on the curve $C_f(\mathbb{F}_p)$. Suppose that the first two columns of $A$ are linearly independent. Then, the cardinality of the set
 $ \{ (k,l) \in \mathbb{F}_p^2 : \nu_p(f(r_1+kp, s_1+lp)) =1 \text{ and }\nu_p(f(r_2+kp, s_2+lp)) =1\}$ is $(p-1)^2 $, the cardinality of the set
 $ \{ (k,l) \in \mathbb{F}_p^2 : \nu_p(f(r_1+kp, s_1+lp)) >1 \text{ and }\nu_p(f(r_2+kp, s_2+lp)) =1\}$ is $p-1$, and the cardinality of the set
  $ \{ (k,l) \in \mathbb{F}_p^2 : \nu_p(f(r_1+kp, s_1+lp)) >1 \text{ and }\nu_p(f(r_2+kp, s_2+lp)) >1\}$ is $1$.
 \end{proposition}
\begin{proof}
 For  $i,j \in \{0,1,2,\ldots,p-1\}$, consider the following Taylor series expansion,
\begin{equation}\label{11}
	f(r_1+ip,s_1+jp) \equiv f(r_1,s_1)+ipf_x(r_1,s_1)+jpf_y(r_1,s_1) \pmod{p^2}. 
\end{equation} 
We know that $f(r_1,s_1) \equiv 0 \pmod p$ which implies that there exist integers $\alpha$ and $\alpha'$ such that $f(r_1,s_1) = \alpha p$ and $f(r_1+ip,s_1+jp) = \alpha' p$.  So, by using equation \eqref{11}, we have \begin{equation}
	i f_x(r_1,s_1)+j f_y(r_1,s_1) \equiv -\alpha + \alpha'\pmod p. \label{12}
\end{equation} 
Similarly, by using $\nu_p(f(r_2,s_2))=1$ we have 
\begin{equation}
	if_x(r_2,s_2)+jf_y(r_2,s_2) \equiv -\beta + \beta'\pmod p, \label{14}
\end{equation}
for integers $\beta$ and $\beta'$ such that $f(r_2,s_2) = \beta p$ and $f(r_2+ip,s_2+jp) = \beta' p$. 
From equations \eqref{12} and \eqref{14} we get the following system of two equations in variables $i$ and $j$:
\begin{equation}\label{15}
	A \begin{pmatrix}
		i  \\
		j \\
		1
	\end{pmatrix}\equiv\begin{pmatrix}
		\alpha'  \\
		\beta' 
	\end{pmatrix}\pmod p.\end{equation}
 Since the first two columns of $A$ are linearly independent, the above system of equations has a unique solution for each $\alpha'$ and $\beta'$.  
 \par Now the set $ \{ (k,l) \in \mathbb{F}_p^2 : \nu_p(f(r_1+kp, s_1+lp)) =1 \text{ and }\nu_p(f(r_2+kp, s_2+lp)) =1\}$ corresponds to $\alpha' \not \equiv 0 \pmod p$ and $\beta' \not \equiv 0 \pmod p$. Hence, its cardinality is $(p-1)^2$.  Similarly, the set $ \{ (k,l) \in \mathbb{F}_p^2 : \nu_p(f(r_1+kp, s_1+lp)) >1 \text{ and }\nu_p(f(r_2+kp, s_2+lp)) =1\}$ corresponds to $\alpha'  \equiv 0 \pmod p$ and $\beta' \not \equiv 0 \pmod p$. So, this set has $p-1$ elements. Finally, the set $ \{ (k,l) \in \mathbb{F}_p^2 : \nu_p(f(r_1+kp, s_1+lp)) >1 \text{ and }\nu_p(f(r_2+kp, s_2+lp)) >1\}$ corresponds to $\alpha'  \equiv 0 \pmod p$ and $\beta' \equiv 0 \pmod p$. So, it has only one element.
 \end{proof}
 
 Next, in the case where the first two columns of $A$ are linearly dependent, we establish the following.
\begin{proposition}\label{prop4}
 Let $p$ be a prime number and let $f(x,y)\in\mathbb{Z}[x,y]$ be a polynomial with no solutions to the set of equations $f(x,y)\equiv 0 \pmod {p^2}, f_x(x,y)\equiv 0 \pmod p \text{ and } f_y(x,y)\equiv 0 \pmod p$. Let $(r_1,s_1)$ and $(r_2,s_2)$ be points on the curve $C_f(\mathbb{F}_p)$. Suppose that the first two columns of $A$ are linearly dependent and at least one of them is non-zero. Then, the cardinality of the set
$\{ (k,l) \in \mathbb{F}_p^2 : \nu_p(f(r_1+kp, s_1+lp)) =1 \text{ and }\nu_p(f(r_2+kp, s_2+lp)) =1\}$ is $p(p-1)$ if $r=1$ and $p(p-2)$ if $r=2$, the cardinality of the set $\{ (k,l) \in \mathbb{F}_p^2 : \nu_p(f(r_1+kp, s_1+lp))>1 \text{ and }\nu_p(f(r_2+kp, s_2+lp)) =1\}$ is $0$ if $r=1$ and $p$ if $r=2$, and the cardinality of the set $\{ (k,l) \in \mathbb{F}_p^2 : \nu_p(f(r_1+kp, s_1+lp)) >1 \text{ and }\nu_p(f(r_2+kp, s_2+lp)) >1\}$ is $p$ if $r=1$ and $0$ if $r=2$.
\end{proposition}
\begin{proof}
We begin the proof by noting that equation \eqref{15} can be established just as in Proposition \ref{prop3}.
We also note that $r\neq 0$ by the given assumption on $f(x,y)$. 

\textit{Case 1.} The set $\{ (i,j) \in \mathbb{F}_p^2 : \nu_p(f(r_1+ip, s_1+jp)) =1 \text{ and }\nu_p(f(r_2+ip, s_2+jp)) =1\}$ corresponds to $\alpha' \not \equiv 0 \pmod p$ and $\beta' \not \equiv 0 \pmod p$.  
If $r=1$, then the system \eqref{15} has solutions if and only if $\begin{pmatrix}
	\alpha'  \\
	\beta' 
\end{pmatrix}$ is a multiple of any non-zero column of $A$ and this has $p-1$ choices. For each such choice, 
there are $p$ solutions of \eqref{15}, as the first two columns of $A$ are linearly dependent. Hence, the system \eqref{15} has $p(p-1)$ solutions.

If $r=2$, then 
the system \eqref{15} has solutions if 
and only if $\begin{pmatrix}
	\alpha' -f(r_1,s_1)/p \\
	\beta' - f(r_2,s_2)/p 
\end{pmatrix}$ is a multiple of any of the first two (non-zero) columns of $A$. This has $p-2$ choices, as there are two distinct multiples which lead to $\alpha'=0$, $\beta'=0$ respectively.
Hence, the system \eqref{15} has $p(p-2)$ solutions.

\textit{Case 2.} The set $\{ (i,j) \in \mathbb{F}_p^2 : \nu_p(f(r_1+ip, s_1+jp)) >1 \text{ and }\nu_p(f(r_2+ip, s_2+jp)) =1\}$ corresponds to $\alpha'  \equiv 0 \pmod p$ and $\beta' \not\equiv 0 \pmod p$. So, the system \eqref{15} becomes \begin{equation}\label{17}
	A \begin{pmatrix}
		i  \\
		j \\
		1
	\end{pmatrix}\equiv\begin{pmatrix}
		0  \\
		\beta' 
	\end{pmatrix}\pmod p.\end{equation}
Note that if $r=1$, the system \ref{17} has solutions only if the first row of $A$ has all zero entries, which is ruled out by the assumption on $f(x,y)$. Hence, the system \ref{17} has zero solutions when the rank of $A$ is $1$ and $p$ solutions when the rank of $A$ is $2$.

\textit{Case 3.} The set $\{ (i,j) \in \mathbb{F}_p^2 : \nu_p(f(r_1+ip, s_1+jp)) >1 \text{ and }\nu_p(f(r_2+ip, s_2+jp)) >1\}$ corresponds to $\alpha'  \equiv 0 \pmod p$ and $\beta' \equiv 0 \pmod p$. So, the system \eqref{15} becomes \begin{equation}\label{16}
	A \begin{pmatrix}
		i  \\
		j \\
		1
	\end{pmatrix}\equiv\begin{pmatrix}
		0  \\
		0 
	\end{pmatrix}\pmod p.\end{equation}
Hence, the system \ref{16} has $p$ solutions when $r=1$, and zero solutions when $r=2$.
\end{proof}

\section{Indicator random variables for valuation greater than $1$}\label{sec3}

We shall assume for the rest of this work that $p$ is a prime number and that $f$ satisfies the condition in Theorem \ref{thm1}, namely that there are no solutions to the set of equations $f(x,y)\equiv 0 \pmod p, f_x(x,y)\equiv 0 \pmod p,$ and $ f_y(x,y)\equiv 0 \pmod p$ for all sufficiently large prime numbers $p$.

For each $(x,y) \in C_f(\mathbb{F}_p)$, define a random variable $Y_{x,y} : \mathbb{F}_p \times \mathbb{F}_p \rightarrow \{0,1\}$ as
        \[
   Y_{x,y}(k,l)= 
\begin{cases}
   1,& \nu_p(f(x+kp,y+lp)) >1, \hspace*{.2cm}  \\
    0,              & \text{otherwise}.
\end{cases}
\]
 In other words, $Y_{x,y}$ is an indicator function on the translates of $(x,y)$ where $f$ has $p$-adic valuation greater than 1. 
  For all points $(x,y) \in C_f(\mathbb{F}_p)$, at least one of $f_x(x,y)$ or $ f_y(x,y)$ is non-zero in $\mathbb{F}_p$ by the assumption on $f$. Hence, using Proposition \ref{Pre 1}, we infer that the probability mass function of $Y_{x,y}$ is
  \begin{equation}\label{pmf}
  	P(Y_{x,y}=i)=
  \begin{cases}
	1/p,  & i=1,  \\
	(p-1)/p,  & i=0.
  \end{cases}
  \end{equation}
  Note that $Y_{x,y}$ has the following properties: 
   \begin{enumerate}
   \item[(i)] The expectation $E[Y_{x,y}]=1/p.$ 
   \item[(ii)] The variance $\operatorname{Var}(Y_{x,y}) = \frac{p-1}{p^2}$.
   \item[(iii)] The $r^{\mathrm{th}}$ moment about the origin $\mu_r' = 1/p$.
   \item[(iv)] By Chebyshev’s inequality 
   	$P(|Y_{x,y}-1/p| \geq \epsilon) \leq \dfrac{p-1}{\epsilon^2 p^2},$
    for any $\epsilon>0$.
   \end{enumerate}
 
 For any pair of points $(x,y)$ and $(r,s) \in C_f(\mathbb{F}_p)$ that satisfy the conditions of Proposition \ref{prop3}, the joint probability mass function of $Y_{x,y}, Y_{r,s}$ is \[
   P(Y_{x,y}=i,Y_{r,s}=j)= 
\begin{cases}
   1/p^2,& i=j=1, \\
    (p-1)/p^2, & i=1,j=0 \text{ or } i=0,j=1, \\
    (p-1)^2/p^2, & i=j=0. 
\end{cases}
\]  
Clearly $P(Y_{x,y}=i,Y_{r,s}=j) = P(Y_{x,y}=i)P(Y_{r,s}=j)$ for all $i$ and $j$, as \eqref{pmf} holds for both points.  Hence, $Y_{x,y}$ and $Y_{r,s}$ are independent random variables. 

Similarly, the joint probability mass function for $Y_{x,y}, Y_{r,s}$ where the pair of points satisfy the conditions of Proposition \ref{prop4} is
\[
   P(Y_{x,y}=i,Y_{r,s}=j)= 
\begin{cases}
	0,& i=j=1, \\
	1/p, & i=1,j=0 \text{ or } i=0,j=1, \\
	(p^2-2 p)/p^2, & i=j=0. 
\end{cases}
\]
in the case $r=2$. In the other case $r=1$, we have
\[
   P(Y_{x,y}=i,Y_{r,s}=j)= 
\begin{cases}
	1/p,& i=j=1, \\
	0, & i=1,j=0 \text{ or } i=0,j=1, \\
	(p^2- p)/p^2, & i=j=0. 
\end{cases}
\]
Clearly, neither of these cases describe independent random variables. However, the correlation coefficient is $-1/(p-1)$ for the case $r =2$ and $1$ for the case $r=1$. Hence, the random variables are close to being uncorrelated when $r=2$, but display a linear relationship when $r=1$.
 
 For ease of notation, we now fix any order on $C_f(\mathbb{F}_p)$ for the rest of this article and denote by $Y_i$ the random variable corresponding to the $i^{\mathrm{th}}$ point.

Recall from Section \ref{sec1} that $m_{k}(f,p)$ denotes the number of element-wise distinct $k$-tuples of points on $C_f(\mathbb{F}_p)$ whose associated matrix $B_k(\mathbf{P})$, given by \eqref{B_k} has rank $2$ with the first two columns of $B_k(\mathbf{P})$ being linearly independent. Let $o_k(f,p)$ denote the number of element-wise distinct $k$-tuples where $B_k(\mathbf{P})$ has rank $1$ with at least one among the first two columns non-zero.
\begin{proposition}\label{prop5} For $k\geq 1$, we have
	\begin{equation}\label{eq-prop5}
		\sum_{l_1 \neq l_2 \neq \cdots \neq l_k} E[Y_{l_1}Y_{l_2} \cdots Y_{l_k}] = m_{k}(f,p)/p^2 + o_k(f,p)/p.
	\end{equation}
\end{proposition}
\begin{proof}
	For any element-wise distinct tuple $l_1, l_2,\ldots, l_k$, it is easy to see that $P(Y_{l_1}=Y_{l_2}=\cdots=Y_{l_k}=1)=n/p^2$ where $n$ is the number of solutions $(i,j)$ of the equation 
	 \begin{equation}\label{no-of-solns}
		B_k(\mathbf{P}) \begin{pmatrix}
			i  \\
			j \\
			1
		\end{pmatrix}\equiv\begin{pmatrix}
			0  \\
			0 \\
			\vdots \\
			0
		\end{pmatrix}\pmod p.\end{equation}
In the case where the rank of $B_k(\mathbf{P})$ is $2$ and the first two columns are linearly independent, there is a unique solution of \eqref{no-of-solns} and hence $n=1$. In the case where the rank of $B_k(\mathbf{P})$ is $1$ with at least one among the first two columns non-zero, there are $n=p$ distinct solutions of \eqref{no-of-solns}. The case where the rank of $B_k(\mathbf{P})$ is $0$ is ruled out by the condition on $f$ in Theorem \ref{thm1}. 

The remaining cases are those where the rank of $B_k(\mathbf{P})$ is either $3$, or $2$ with the first two columns linearly dependent, or $1$ with the first two columns zero. It can be seen that there are no solutions of \eqref{no-of-solns} in either of these remaining cases.

Since $E[Y_{l_1}Y_{l_2} \cdots Y_{l_k}]= P(Y_{l_1}=Y_{l_2}=\cdots=Y_{l_k}=1)$, equation \eqref{eq-prop5} follows.
\end{proof}

\begin{corollary}\label{cor1}
	For $k\geq 1$, assuming that Conjecture \ref{conj1} holds, we have
	\begin{align*}
		\lim_{p\rightarrow\infty}\sum_{l_1 \neq l_2 \neq \cdots \neq l_k} E[Y_{l_1}Y_{l_2} \cdots Y_{l_k}] = 1
	\end{align*}
\end{corollary}
\begin{proof}
Recall that $m_{1}(f,p)=0$ trivially and $o_{1}(f,p)=m$. Hence, by the Hasse-Weil bounds \cite{eight}, we get $m_{1}(f,p)/p^2 \rightarrow 0$ and $ o_{1}(f,p)/p \rightarrow 1$ as $p\rightarrow\infty$. Using this in equation \eqref{eq-prop5}, we get the required result for $k=1$.

Also, for $k\geq 2$, we note that $o_k(f,p)$ is bounded by a constant depending on $f$. This is because, for $k=2$, we can bound $o_{k}(f,p)$ by the number of solutions of the following system of $4$ equations in $4$ unknowns
\begin{align*}
	f(x,y) \equiv f(r,s) &\equiv 0 \pmod p; \\
	f_x(x,y)f_y(r,s)-f_x(r,s)f_y(x,y) &\equiv 0 \pmod p; \\
	f_x(x,y)f(r,s)-f_x(r,s)f(x,y) &\equiv 0  \pmod{p^2},
\end{align*}
 and Bezout's theorem gives a bound on this number depending only on $f$ and independent of $p$. For the case $k>2$, there are more equations than there are unknowns and, in general, there are no solutions. In either case, $o_k(f,p)/p\rightarrow 0$ as $p\rightarrow\infty$. Using this and Conjecture \ref{conj1} in equation \eqref{eq-prop5}, we get the required result for $k\geq 2$.
\end{proof}

\section{Uniform distribution results}\label{sec4}
 We start by describing the idea of uniform distribution of the sequence $\mathcal{F}_p$ (defined in \eqref{Fcal_defn}) in $\mathbb{F}_p^3$, namely that for every set $E\subseteq \mathbb{F}_p^3$, the proportion of points of $\mathcal{F}_p$ that lie in $E$ is approximated by the proportion of $|E|$ to $|\mathbb{F}_p^3|$ for $p$ sufficiently large, viz.
 \begin{equation}\label{unif}
    \frac{|\mathcal{F}_p \cap E|}{m}  \approx \frac{|E|}{p^3}, \quad p\rightarrow\infty.
 \end{equation}
 We shall use several measures to quantify the deviation from uniform distribution, also called discrepancy in the literature (see \cite{drnota} for an excellent exposition) and define precise versions of uniform distribution.
 
 At the outset, we describe the version of uniform distribution that provably implies Conjecture \ref{conj1}. For that, we give the following definition of discrepancy.
\begin{equation}\label{D_p}
	\mathcal{D}_p(\mathcal{F}_p) = \sup_{E\subseteq \mathbb{F}_p^3} \left|\frac{|\mathcal{F}_p \cap E|}{|E|}\frac{p^3}{m} - 1 \right|.
\end{equation}
Using this, we formulate the following conjecture, which implies the uniform distribution of $\mathcal{F}_p$ as described in \eqref{unif}.
\begin{conjecture}\label{conj2}
	For $p\rightarrow \infty$, we have $\mathcal{D}_p(\mathcal{F}_p)\rightarrow 0$.
\end{conjecture}
We now give the proof of Conjecture \ref{conj1} assuming this version of uniform distribution.
\begin{proposition}
	 Conjecture \ref{conj2} implies Conjecture \ref{conj1}.
\end{proposition}
\begin{proof}
We start by noting that Conjecture \ref{conj2} can be applied to any independent choice of $k$ vectors from $\mathcal{F}_p$. More precisely, Conjecture \ref{conj2} implies
\begin{equation}\label{disc}
	\sup_{E\subseteq (\mathbb{F}_p^3)^k} \left|\frac{|\mathcal{F}_p^{k'} \cap E|}{|E|}  \frac{p^{3k}}{m(m-1)\ldots(m-k+1)}-1 \right| \rightarrow 0,
\end{equation}	
as $p\rightarrow\infty$ where $\mathcal{F}_p^{k'}$ denotes the element-wise distinct $k$-tuples of elements of $\mathcal{F}_p$.

Applying \eqref{disc} to $E=\{(v_1,v_2,\ldots,v_k) \in (\mathbb{F}_p^3)^{k} |\ \textit{associated matrix has rank } 2$ \textit{ and first two columns linearly independent}$\}$ gives
\begin{equation*}
	\left|\frac{m_{k}(f,p)}{|E|}\frac{p^{3k}}{m(m-1)\cdots(m-k+1)} - 1\right| \rightarrow 0.
\end{equation*} 
Now $|E| = p^3(p^3-p)(p^2-2)(p^2-3)\cdots(p^2-k+1)$ and $m/p\rightarrow 1$ as $p\rightarrow\infty$ by the Hasse-Weil bounds \cite{eight}. Hence, we get (as $p\rightarrow\infty$)
\begin{equation*}\label{claim1}
	\left|\frac{m_{k}(f,p)}{p^2} - 1\right| \rightarrow 0,
\end{equation*} 
which is the same as \eqref{lim1}.
\end{proof}

 For comparison with Conjecture \ref{conj2}, we recall the standard version of discrepancy in the $3$-dimensional unit cube of the sequence of $\mathcal{F}_p$ scaled down by a factor of $p$ in each component, which we denote $\Delta_p(\mathcal{F}_p)$.
\begin{equation}\label{Delta_p}
	\Delta_p(\mathcal{F}_p) = \sup_{E\subseteq \mathbb{F}_p^3} \left|\frac{|\mathcal{F}_p \cap E|}{m} - \frac{|E|}{p^3}\right|.
\end{equation}

We now state the following conjecture, which is another version of uniform distribution as described in \eqref{unif}.
\begin{conjecture}\label{conj-new}
	For $p\rightarrow\infty$, we have $\Delta_p(\mathcal{F}_p)\rightarrow 0$.
\end{conjecture}
 It can be seen that Conjecture \ref{conj2} implies Conjecture \ref{conj-new}, as $\Delta_p(\mathcal{F}_p) \leq \mathcal{D}_p(\mathcal{F}_p)$. It is also easy to see that Conjecture \ref{conj-new} gives $m_{2,2}/p^2 \rightarrow 1$ as $p\rightarrow\infty$. However, it is not yet clear whether this suffices to prove Conjecture \ref{conj1} in its entirety. 
 
 On the other hand, we note the connection of Conjecture \ref{conj-new} with exponential sums by recalling the Erd\"os-Turan-Koksma inequality (Theorem 1.21, \cite{drnota}).
\begin{proposition}\cite{drnota}\label{lem1}
There exists a constant $C>0$ such that, for any positive integer $L$, the following bound holds
\begin{equation*}
	\Delta_p(\mathcal{F}_p) \leq C\left(\frac{1}{L} + \frac{1}{m}\sum_{0<|\mathbf{a}|\leq L} \frac{1}{r(\mathbf{a})}
	\left|\sum_{\mathbf{v}\in\mathcal{F}_p}e_p(\mathbf{a}\cdot \mathbf{v})\right|
	\right) ,
\end{equation*}
where the sum runs over triples $\mathbf{a}\in \mathbb{Z}^3$ which satisfy the given condition with $|\mathbf{a}|=\max_{i=1}^3 a_i$, $r(\mathbf{a})=\prod_{i=1}^3\max(|a_i|,1)$ and $e_p(x) = \exp(2\pi i x/p)$.
\end{proposition}
We finally state the following conjecture on exponential sums which appears to be obtainable along the same lines as Bombieri \cite{bomb}.
\begin{conjecture}\label{conj3}
There exists a constant $D>0$ such that, for any prime $p$ and any $k,l\in \{0,1,\dots,p-1\}$,  the following bound holds
\begin{equation*}
	\sum_{(x,y)\in C_f(\mathbb{F}_p)} e_{p^2}\left(f(x+kp,y+lp)\right) \leq D p^{\frac{1}{2}}.
\end{equation*}
\end{conjecture}
It is easy to see, using Proposition \ref{lem1}, that Conjecture \ref{conj3} implies $\Delta_p(\mathcal{F}_p) \ll p^{-\frac{1}{2}+\epsilon}$  for any $\epsilon>0$, which definitely implies Conjecture \ref{conj-new}. 

Numerical computations seem to support the conjectures on uniform distribution  but the speed of convergence appears to follow Conjecture \ref{conj3}. More precisely, we compute lower bounds for $\Delta_p(\mathcal{F}_p)$ and $ \mathcal{D}_p(\mathcal{F}_p)$ in Table \ref{tab:tab1} that are comparable to the order of decay predicted by Conjecture \ref{conj3}.  These are obtained by restricting $E$ to cubes of side $p/3$ in \eqref{D_p} and \eqref{Delta_p}.
\renewcommand{\arraystretch}{1.9}
\begin{table}[h!]
\centering
\begin{tabular}{|c|c|c|c|}
\hline
$p$ & Lower bound for $\Delta_p$ & Lower bound for $\mathcal{D}_p$ & $p^{-1/2}$ \\
\hline
$1009 $ & $0.01088009538$ & $0.29376257545$& $0.03148142750$\\
\hline
$2003$ & $0.01131027426$ & $0.30537740503$& $0.02234392810$ \\
\hline
$3001$ & $0.00930693824$ & $0.25128733264$& $0.01825437644$\\
\hline
$4003$ & $0.00494879436$ & $0.13361744798$& $0.01580546236$ \\
\hline
$5443$ & $0.00511927708$ & $0.13822048125$& $0.01355441669$ \\
\hline
$6007$ & $0.00672228668$ & $0.18150174042$& $0.01290242026$ \\
\hline
$7013$ & $0.00646586788$ & $0.17457843276$& $0.01194120297$ \\
\hline
$8011$ & $0.00496738448$ & $0.13411938098$& $0.01117266132$ \\
\hline
$9007$ & $0.00563262285$ & $0.15208081695$& $0.01053682867$ \\
\hline
$10009$ & $0.00494777953$ &$0.13359004739$ & $0.00999550303$ \\
\hline
\end{tabular}
\caption{Lower bounds on discrepancy for various primes for the polynomial $f(x,y)=x^3+y^2+x y +1$.}
\label{tab:tab1}
\end{table}

\section{Proof of Theorem \ref{thm1}}\label{sec5}

We start by noting that the random variable $X(k,l)$ defined in \eqref{Xdefn}, representing the total number of points in $C_f(\mathbb{F}_p)$ with valuation greater than $1$ in a random $p \times p$ square, has the following expression.
\begin{equation*}\label{X_def}
	X(k,l)= \sum_{(x,y)\in C_f(\mathbb{F}_p)}Y_{x,y}(k,l)=\sum_{i=1}^m Y_i(k,l).
\end{equation*}
Recall that $m=|C_f(\mathbb{F}_p)|$. 

By the properties of the random variables $Y_{x,y}$ discussed in Section \ref{sec3}, we have 
\begin{align*}
E[X] &= \sum_{i=1}^{m}E[Y_i] =  m/p. 
\end{align*} 

In general for $k\geq 1$
\begin{align}\label{mom_k}
	E[X^k] &=E\big[\big(\sum_{l=1}^m Y_{l}\big)^k\big] = E\big[\big(\sum_{l_1,l_2,\cdots,l_k=1}^m Y_{l_1}Y_{l_2} \ldots Y_{l_k}\big)\big]  \notag\\
	&= \sum_{i=1}^k S(k,i)\sum_{l_1 \neq l_2 \neq \cdots \neq l_i} E[Y_{l_1}Y_{l_2} \ldots Y_{l_i}],
\end{align}
where $l_i \in \{1,2,\dots,m\}$ and $S(k,i)$ is the Stirling number of the second kind. This is because for each element-wise distinct tuple $(l_1,l_2,\ldots,l_i)$ in the last sum, the contribution from each product, where only the variables $Y_{l_1},Y_{l_2},\ldots,Y_{l_i}$ appear, is the value  $E[Y_{l_1}Y_{l_2} \ldots Y_{l_i}]$. Since the positions where these variables appear, form a partition of $\{1,2,\cdots,k\}$ into $i$ parts, we get $S(k,i)$ many such products for each $(l_1,l_2,\ldots,l_i)$.

 Thus, from \eqref{mom_k} and Corollary \ref{cor1} (assuming Conjecture \ref{conj1} as in the theorem), we get the following limit for the $k^{\mathrm{th}}$ moment $ E[X^k]$: 
\begin{equation*}
\lim_{p\to\infty} E[X^k] = \sum_{i=1}^k S(k,i),
\end{equation*}
which is indeed the $k^{\mathrm{th}}$ moment of the Poisson distribution with parameter $\lambda=1$ (see \cite{six}). In other words, the moment generating function of $X$, $M_X(t)$, converges to
\begin{equation*}
 	\lim_{p\to\infty}  M_X(t) = \exp(e^t - 1).
\end{equation*}
Hence, $X$ converges in distribution to a Poisson random variable with parameter $\lambda=1$ as $p\rightarrow \infty$.
\section*{Acknowledgements}
The authors would like to thank the anonymous referee for valuable suggestions and
insightful comments that  helped to vastly improve the exposition.

\end{document}